\documentclass[draft]{article}

\usepackage{amsmath}
\usepackage{amsthm}
\usepackage{amscd}
\usepackage{amssymb}
\usepackage{latexsym}

\newcommand{\LL}{{\cal L}}
\newcommand{\CA}{{\cal A}}
\newcommand{\CB}{{\cal B}}

\newcommand{\LG}{\mathfrak L}
\newcommand{\DB}{\|}
\newcommand{\wbar}{\overline}

\newcommand{\chop}{\dagger}

\newcommand{\dvb}{\|}
\def\epsilon{\varepsilon}

\newcommand{\Out}{\mbox{Out}}
\newcommand{\Aut}{\mbox{Aut}}
\newcommand{\Stab}{\mbox{Stab}}
\newcommand{\inv}{^{-1}}
\newcommand{\BBT}{\mbox{BBT}}
\newcommand{\Ax}{\mbox{Ax}}

\newcommand{\vol}{\mbox{vol}}

\newcommand{\FN}{F_N}   
\newcommand{\cvn}{\mbox{cv}_N}
\newcommand{\CVN}{\mbox{CV}_N}
\newcommand{\cvnbar}{\overline{\mbox{cv}}_N}
\newcommand{\CVNbar}{\overline{\mbox{CV}}_N}
\newcommand{\CQ}{{\cal Q}}

\newcommand{\R}{\mathbb R}
\newcommand{\Z}{\mathbb Z}

\newcommand{\N}{\mathbb N}
\newcommand{\Hy}{\mathbb H}
\newcommand{\A}{{\cal A}}

\def\QED{\hfill\rlap{$\sqcup$}$\sqcap$\par\bigskip}
\def\qed{\hfill\rlap{$\sqcup$}$\sqcap$\par}
\def\bar{\overline}
\def\tilde{\widetilde}

\newtheorem{thm}{Theorem}[section]
\newtheorem{cor}[thm]{Corollary}
\newtheorem{lem}[thm]{Lemma}
\newtheorem{prop}[thm]{Proposition}

\theoremstyle{definition}
\newtheorem{defn}[thm]{Definition}
\newtheorem{rem}[thm]{Remark}

\theoremstyle{remark}

\newtheorem*{acknowledgements}{Acknowledgements}

\newenvironment{proof*}[1][Proof]{}{\QED}

\newcommand{\affiliationone}{}
\newcommand{\email}{\\}
\renewcommand{\and}{ and }

\newcommand{\singlebox}{}
\newcommand{\esinglebox}{}

\title{$\R$-trees and laminations for free groups II: \\ The dual lamination of an $\R$-tree}

\author{Thierry Coulbois, Arnaud Hilion \and Martin Lustig}

\begin{document}

\maketitle

\section{Introduction}\label{sec:intro}

A geodesic lamination $\LG$ on a hyperbolic surface $S$, provided with
a transverse measure, defines (via the lift $\tilde \LG$ of $\LG$ to
the universal covering of $S$) an action of $\pi_{1} S$ on an
$\R$-tree $T$ which is often called {\em dual} to the lamination
$\tilde \LG$. Conversely, every small action of $\pi_{1} S$ on an
$\R$-tree $T$ comes from this construction, provided the surface is
closed and the action on $T$ is {\em small} \cite{skora}.

\smallskip

A generalization of this concept occurs first in the work of E.~Rips,
and is since then widely used.  A particular kind of $\R$-tree actions
can be defined as spaces dual to measured foliations on finite
$2$-complexes. If the latter is not a surface (and can not be made
into a surface by certain elementary moves), then the resulting
$\R$-tree is qualitatively different from the ones dual to surface
laminations, see \cite{glp,bf95,lp}.

\smallskip

A third occurence of laminations in direct relation to $\R$-trees
takes place in the context of free group automorphisms, specifically
for {\em irreducible automorphisms with irreducible powers (= iwip
automorphisms)}.  For such an $\alpha \in \Aut(\FN)$ every
(non-periodic) conjugacy class of elements in $\FN$ converges to a
collection of biinfinite {\em legal} paths on a train track
representing $\alpha$, see \cite{luhabil,bfh}.  On the other hand
(compare \cite{bf?}, \cite{gjll}), the train track itself converges
towards an $\R$-tree $T_{\alpha}$ which is projectively fixed under
the induced action of $\alpha$ on the set of $\R$-tree actions of
$\FN$.  In the particular case where $\alpha$ is induced by a
pseudo-Anosov homeomorphism $\varphi$ of some surface $S$ (with one
boundary component), then the above collection of biinfinite legal
paths is in precise 1~-~1 correspondence with the leaves of the stable
lamination of $\varphi$.

\smallskip

Finally, in \cite{ll4} a collection of one-sided infinite words
reminiscent to half-leaves of a lamination was associated to an
arbitrary very small $\FN$-action on an $\R$-tree $T$, as a tool to
prove that iwip automorphisms of $\FN$ have a North-South dynamics on
the space $\CVNbar$ of very small $\R$-tree actions of $\FN$.

\smallskip

This puzzle of coinciding and consistent observations induced the
authors to set out for a general theory, in the realm of free (and
perhaps later word-hyperbolic) groups. As a first step, in
\cite{chl1-I} algebraic laminations were defined, generalizing at the
same time geodesic laminations on surfaces, as well as symbolic flows
as known from discrete dynamics. They come in three equivalent
languages, group theoretic, dynamic and combinatorial, and passing
from one to the other turns out to be rather helpful. In \cite{chl1-I}
these ``translations'' were established with some care, and the
topology, the partial order, as well as the natural $\Out(\FN)$-action
were studied.

\bigskip

In the present paper we use these tools to define, for any isometric
$\FN$-action on an $\R$-tree $T$, a set
\[
L^{2}(T) \subset \partial^{2}\FN := 
\partial \FN \times \partial \FN \smallsetminus 
\Delta \, ,
\] 
where $\partial \FN$ denotes the Gromov boundary of the free group
$\FN$, and $\Delta$ is the diagonal.  The set $L^{2}(T)$ is empty if
the $\FN$-action on $T$ is free and discrete (simplicial), and it is
an algebraic lamination otherwise.  There are several competing
natural approaches to define $L^2(T)$, which we present briefly below.
Working out the precise relationship between them is the core content
of this paper.

\medskip

\noindent{\bf 1.}  (see \S\ref{subsec:lamtree}) The lamination $L^{2}_\Omega(T)$ is
defined for all isometric actions of $\FN$ on an $\R$-tree $T$: For
every $\epsilon > 0$ we consider the set $\Omega_{\epsilon}(T)$ of all
elements $g \in \FN$ with translation length on $T$ that satisfies:
\[
\dvb g \dvb_{T} \, \, < \, \epsilon
\]
The set $\Omega_{\epsilon}(T)$ generates an algebraic lamination
$L^{2}_{\epsilon}(T)$ which is the smallest lamination that contains
every $(g\inv g\inv g\inv \ldots, ggg\ldots) \in \partial^{2}\FN$ with
$g \in \Omega_{\epsilon}(T)$. We define $L_{\Omega}^{2}(T)$ to be the
intersection of all $L^{2}_{\epsilon}(T)$.

\medskip

We are most interested in $\R$-trees $T$ where every $\FN$-orbit of
points is dense in $T$.  To any such $T$ we associate in this paper
two more laminations, which are of rather different nature:

\medskip

\noindent{\bf 2.} (see \S\ref{subsec:lamtreeonesided}) In order to
define the lamination $L^2_{\infty}(T)$ one first fixes a basis $\CA$
of $\FN$. Then one picks an arbitrary point $P \in T$ and considers
the set $L_{\CA}^{1}(T)$ of one-sided infinite reduced words $x_{1}
x_{2} \ldots$ in the basis $\CA$ such that the set of all $x_{1} x_{2}
\ldots x_{k} P$ has bounded diameter in $T$.  One immediately observes
that the set $L_{\CA}^{1}(T)$ is independent of the choice of $P$.  As
next step, one considers the language $\LL_{\CA}(L_{\CA}^{1}(T))$
derived from $L_{\CA}^{1}(T)$, i.e. the set of all finite subwords of
any $x_{1} x_{2} \ldots \in L_{\CA}^{1}(T)$, and its {\em recurrent}
sublanguage $\LL^\infty_{\CA}(L_{\CA}^{1}(T))$: The latter consists
precisely of those words which occur infinitely often as subword in
some $x_{1} x_{2} \ldots \in L_{\CA}^{1}(T)$.  The advantage of this
language is that it is {\em laminary}, and thus it defines canonically
an algebraic lamination $L^{2}_{\infty}(T)$ (see \cite{chl1-I}).  We
prove in detail in \S\ref{subsec:lamtreeonesided} that this algebraic
lamination does not depend on the basis $\CA$ used in the construction
sketched above.

\medskip

\noindent{\em Aside:} The passage from $L^{1}(T)$ to the recurrent
language $\LL^\infty_{\CA}(L_{\CA}^{1}(T))$ is subtle but rather
important: It turns out that, contrary to $L^{2}_{\infty}(T)$, the set
$L^{1}(T)$ depends not just on the topology, but actually on the
metric of $T$. This is discussed in detail in \cite{chl2} and
\cite{chll}.

\medskip

\noindent{\bf 3.} (see \S\ref{sec:LQetQ2}) If $T$ is an $\R$-tree dual
to a measured lamination $\mathfrak L$ on surface $S$, then a leaf $l$
of the lift $\tilde {\mathfrak L}$ of $\mathfrak L$ to the universal
covering $\tilde S$ determines on one hand a point $x_{l}$ in the dual
tree $T$, and on the other hand two limit points $P_{+}, P_{-} \in
\partial \tilde S$ on the boundary at infinity $\partial \tilde S$ of
$\tilde S$:
\[
\{P_{+}, P_{-}\} = \partial l\, .
\]
Note that, in the case where $S$ has non-empty boundary, $\partial
\tilde S$ is canonically identified with the Gromov boundary $\partial
\FN$ of $\FN = \pi_{1} S$.

In \cite{ll4} this correspondence has been generalized to an
$\FN$-equivariant map $\CQ: \partial \FN \to \bar T \cup \partial T$,
for any $\R$-tree $T$ with dense $\FN$-orbits, where $\bar T$ denotes
the metric completion of $T$, and $\partial T$ the Gromov boundary.
The definition of the map $\CQ$ is reviewed in \S\ref{subsec:bbtq},
and the geometric meaning of $\CQ$ is explained in more detail in
\S\ref{subsec:lamonsur}.  In the above special case one gets:
\[
\CQ(P_{+}) = \CQ(P_{-}) = x_{l} \in T
\]
This motivates the definition of the lamination $L^{2}_{\CQ}(T)$,
which consists of all pairs $(X, X') \in \partial^{2}\FN$ that
determine the same limit point $\CQ(X) = \CQ(X')$ in $\bar T$.

\medskip

The main result of this paper, proved in two steps (Propositions
\ref{prop:langequi} and \ref{prop:lamtree}), is:

\begin{thm}
\label{theoremtwo}
For every very small $\R$-tree with $\FN$-action that has dense orbits
the above described three algebraic laminations coincide:
\[
L^2_\Omega(T) = L^2_\infty(T) = L^2_{\CQ}(T)
\]
This defines a dual algebraic lamination $L^{2}(T)$ canonically
associated to $T$.
\end{thm}

If the $\FN$-action on $T$ does not have dense orbits, we define
$L^{2}(T) = L^{2}_{\Omega}(T)$. For any non-simplicial such $T$ there
is a canonical (maximal) quotient $\R$-tree $T'$ which has dense
orbits, and we show that $L^{2}(T) \subset L^{2}(T')$ (see
Remark~\ref{discrete-dense} and Remark~\ref{contractions}).

In a subsequent third paper we go one step further and consider
invariant measures $\mu$ (called {\em currents}) on $L^{2}(T)$.  We
show in \cite{chl1-III} that such a current defines a {\em dual
metric} $d_{\mu}$ on $T$: If $T$ is dual to a surface lamination $\LG$
as in the beginning of this introduction, and if $\mu$ comes from a
transverse measure on $\LG$, then $d_{\mu}$ is indeed the $\R$-tree
metric on the dual tree $T$. In general, however, it is shown in
\cite{chl1-III} that this dual metric can have very exotic properties.

\smallskip

Note also that the lamination $L^{2}(T)$ introduced in this paper has
been used successfully in \cite{chl2} to characterize the underlying
topological structure of $\R$-trees which stays invariant when the
metric is $\FN$-equivariantly changed (so called ``non-uniquely
ergometric $\R$-trees'').

\smallskip

The dual lamination $L^{2}(T)$ plays also a crucial role in
establishing a continuous and $\Out(\FN)$-equivariant map from a large
part of the boundary of Outer space into the space of projectivized
currents, see \cite{chl3} and \cite{kl1}.  It is also the basis for
work in progress of the third author with I.~Kapovich on perpendicular
($\R$-tree, current)-pairs.

\bigskip

\begin{acknowledgements} This paper originates from a workshop organized
at the CIRM in April 05, and it has greatly benefited from the
discussions started there and continued around the weekly Marseille
seminar ``Teichm\"uller'' (partially supported by the FRUMAM).
\end{acknowledgements}


\section{$\FN$-actions on $\R$-trees}
\label{sec:actions}

In this section we recall some of the known properties of actions of a
free group $\FN$ on an $\R$-tree $T$.  For details and background see
\cite{vogt,shalen} and the references given there.

\smallskip

An $\R$-tree is a metric space which is $0$-hyperbolic and
geodesic. Alternatively, an $\R$-tree is a metric space $(T, d)$ where
any two points $P, Q \in T$ are joined by a unique arc and this arc is
isometric to the interval $[0, d(P, Q)] \subset \R$.

\smallskip

In this paper an $\R$-tree $T$ always comes with a left action of
$\FN$ on $T$ by isometries. Any isometry $w$ of $T$ is either {\em
elliptic}, in which case it fixes at least one point of $T$, or else
it is {\em hyperbolic}, in which case there is an {\em axis} $\Ax(w)$
in $T$, isometric to $\R$, which is $w$-invariant, and along which $w$
acts as translation. The {\em translation length}
\[
\|w\|_T=\inf\{d(P,wP)\,|\,P\in T\}
\]
agrees in the hyperbolic case with $d(Q,wQ)$ for any point $Q$ in
$\Ax(w)$, while in the elliptic case it is $0$.  The $\FN$-action on
$T$ is called {\em abelian} if there exists a homomorphism $\psi: \FN
\to \R$ such that $\|w\|_T = | \psi(w) |$ for all $w \in \FN$.  In
this case there exists an infinite ray $\rho \subset T$ such that
every $w \in \FN$ acts as translation for every point sufficiently far
out on $\rho$.  The $\R$-tree actions treated in this paper will all
turn out to be non-abelian.

\smallskip

We always assume that $T$ is a {\em minimal} $\R$-tree, i.e.  there is
no non-empty $\FN$-invariant proper subtree of $T$.  Another minimal
$\R$-tree $T'$ with isometric $\FN$-action is $\FN$-equivariantly
isometric to $T$ if and only if one has
\[
\|w\|_T= \|w\|_{T'}
\]
for every element $w \in \FN$, and if the actions are non-abelian. The
set of $\R$-trees equipped with such $\FN$-actions inherits a topology
from its image in $\R^{\FN }$ under the map
\[
T \mapsto (\dvb w \dvb_{T})_{w \in \FN} \in \R^{\FN } \, .
\]

\smallskip

A tree (or a tree action) is called {\em small} if any two group
elements that fix pointwise a non-trivial arc in $T$ do commute. It is
called {\em very small} if moreover (i) the fixed set $Fix(g) \subset
T$ of any elliptic element $1 \neq g \in\FN$ is a segment or a single
point (i.e. ``no branching"), and (ii) $Fix(g) = Fix(g^m)$ for all
$g\in\FN$ and $m \geq 1$.  One sees easily that every small (and thus
every very small) action is non-abelian.

\medskip

The particular case of {\em simplicial} $\R$-trees $T$ with isometric
$\FN$-actions, which have trivial edge stabilizers, arises from graphs
$\Gamma$ with a {\em marking isomorphism} $\FN \cong \pi_{1} \Gamma$,
where the edges of $\Gamma$ are given a non-negative length, which is
for at least one of them strictly positive: The simplicial $\R$-tree
$T$ is then given by lifting the edge lengths to the the universal
covering $\tilde \Gamma$, equipped with the action of $\FN$ by deck
transformations.  (Note that there is a minor ambiguity concerning the
the topology on the set $T = \tilde \Gamma$: the metric topology on
$T$ is in general weaker than the cellular topology on $\tilde \Gamma$
!).  If every edge length of $\Gamma$ is non-zero, then the action of
$\FN$ on $T$ is free. The space $\cvn$ of $\R$-trees equipped with a
free simplicial action has been introduced by M.~Culler and
K.~Vogtmann.  Its closure $\cvnbar$ (in the space $\R^{\FN}$ as
described above) is precisely the set of all of the above mentioned
very small $\R$-trees.  The {\em boundary} $\cvnbar \smallsetminus
\cvn$ is denoted by $\partial \cvn$.  One often normalizes $\Gamma$ to
have volume 1, thus obtaining the subspace $\CVN$ of $\cvn$, which has
been named {\em Outer space} by P.~Shalen.  Alternatively, one can
{\em projectivize} the space of tree actions: two trees $T$ and $T'$
are in the same equivalence class $[T]$ if they are
$\FN$-equivariantly homothetic. This projectivization maps $\cvnbar$
onto a compact space $\CVNbar$, which contains a homeomorphic copy of
$\CVN$, called the {\em interior}, and the projectivized image
$\partial\CVN$ of $\partial\cvn$, called the {\em boundary}. Both
$\CVN$ and $\CVNbar$ are contractible and finite dimensional. For more
information see \cite{vogt,CulVog}.

\smallskip

The group $\Out(\FN )$ acts by homeomorphisms on $\cvn$ and
$\partial\cvn$, as well as on $\CVN$ and $\partial \CVN$, and the
action on $\CVN$ is properly discontinuous (though not free). These
actions, specified in \S\ref{subsec:lambdatwo}, provide valuable
information about the group $\Out(\FN )$.  Note also that there is a
strong similarity between Outer space $\CVN$ with the $\Out (\FN
)$-action on one side, and Teichm\"uller space with the action of the
mapping class group on the other. The only substantial difference is
that $\CVN$ is not a manifold.

\medskip

In the second half of this paper we will concentrate on the
particularly interesting case where some (and hence every) $\FN$-orbit
of points is dense in $T$.  That this ``dense orbits'' hypothesis is
not very restrictive follows from the following consideration:


\begin{rem}
\label{discrete-dense}
Every $\R$-tree $T \in \cvn \cup \partial \cvn$ decomposes canonically
into two disjoint $\FN$-invariant (possibly empty) subsets $T_{d}$ and
$T_{c}$, where the former is given as the union of all points $P \in
T$ such that the orbit $\FN P$ is a discrete (or, equivalently, a
closed discrete) subset of $T$, and the latter is the complement
$T\smallsetminus T_{d}$.  Using property (ii) in the above definition
of a very small action, the subset $T_{c} \subset T$ is easily seen to
be closed, and each connected component of it is a subtree $T'$ of $T$
with the property that the subgroup $U$ of $\FN$ that stabilizes $T'$
acts on $T'$ with {\em dense orbits}: $T' = \overline{U P}$ for any $P
\in T'$.  Unless $T_{c}$ is empty, in which case the set of branch
points is a discrete subset of $T$ and thus $T$ is simplicial, we can
contract the closure of every connected component of $T_{d}$ to a
single point, to get the canonical maximal non-trivial quotient tree
$T/T_{d}$ on which now all of $\FN$ acts minimally and with dense
orbits.  Compare \cite{le94}.
\end{rem}


\section{Bounded Back Tracking}
\label{sec:bbt}

Every small action on an $\R$-tree is known to have the {\em bounded
backtracking property (BBT)} (see \cite{gjll}), which is of great use
in this paper:

\smallskip

Let $\Gamma$ be any (non-metric) graph with a marking isomorphism
$\pi_{1} \Gamma\cong \FN$, and let $\tilde \Gamma$ be its universal
covering.  Let $i: \tilde \Gamma \to T$ be any $\FN$-equivariant
map. Then the map $i$ {\em satisfies BBT} if and only if for every
pair of points $P, Q \in \tilde \Gamma$ the $i$-image of the geodesic
segment $[P, Q] \subset \tilde \Gamma$ is contained in the
$C$-neighborhood of $[i(P), i(Q)] \subset T$, where $C \geq 0$ is an a
priori constant independent of the choice of $P$ and $Q$.  We denote
by $\BBT(i) \geq 0$ the smallest such constant.

\smallskip

Every $\R$-tree $T$ with isometric $\FN$-action admits a map $i$ as
above and, $i$ satisfies BBT if and only if any other such map $i':
\tilde \Gamma' \to T$ also satisfies BBT.  Hence the property BBT is a
well defined property of the tree $T$.

We can assume that the above map $i: \tilde \Gamma \to T$ is {\em
edge-geodesic}: $i$ maps every edge $e \subset \tilde \Gamma$ to the
geodesic segment that connects the images of the endpoints of $e$.
One can make $\Gamma$ into a metric graph by giving each edge of
$\Gamma$ and each of its lifts $e$ to $\tilde \Gamma$ the length of
$i(e)$.  Without loss of generality one can assume that the metric on
every edge $e$ is properly distributed so that $i$ is actually {\em
edge-isometric}, i.e. $i$ maps every edge of $\tilde \Gamma$
isometrically onto its image. In this case the inequality
\[
BBT(i) \leq \vol(\Gamma)
\]
has been proved in \cite{gjll}, where the {\em volume}\ $\vol(\Gamma)$
of $\Gamma$ is the sum of the lengths of its edges.

A particular choice of $\Gamma$, for any base
$\cal A$ of $\FN$, is the rose $\cal R_{A}$ with $n$ leaves that are
in
one-to-one correspondence with the elements of $\cal A$. In this case
the universal covering $\cal \tilde R_{A}$ is canonically identified
with the Cayley graph of $\FN$ with respect to the generating
system $\cal A$, and the edge-geodesic map 
$i = i_{P, \cal A}$
is uniquely determined by the choice of a {\em base point} $P = i({\bf
1})$, where ${\bf 1}$ is the vertex of the Cayley graph that
corresponds to the neutral element $1 \in \FN$. In this case we denote
the BBT-constant $\BBT(i_{P,\CA})$ by $\BBT(\CA,P)$, and the volume of
${\cal R}_{\CA}$ by $\vol(\CA, P)$.

\medskip

\begin{lem}\label{lem:PBBT}
Let $T$ be an $\R$-tree with isometric $\FN$-action that satisfies
BBT.  Let $\A$ be any basis of $\FN = F(\CA)$, and let $P$ be any
point of $T$. Then the constants $\BBT(\CA, P) \geq 0$ and
$\mbox{vol}(\CA, P) > 0$ satisfy:

\smallskip

\noindent
(a) For any reduced word $w = x_{1} x_{2} \ldots x_{n}$ in $F(\A)$ and
any prefix $v = x_{1} x_{2} \ldots x_{m}$ of $w$ the point $vP$ is
contained in the $\BBT(\A, P)$-neighborhood of the geodesic segment
$[P, wP] \subset T$.

\smallskip

\noindent
(b) For any cyclically reduced word $w$ in $F(\A)$ one has:
\[
d(wP,P)\leq 2\BBT(\A,P)+\|w\|_T
\]

\smallskip

\noindent
(c) Any subword $u$ of a cyclically reduced word $w \in F(\A)$
satisfies:
\[
d(uP,P)\leq 2\BBT(\A,P)+\|w\|_T
\]

\smallskip

\noindent
(d) Every $x \in \CA \cup \CA\inv$ satisfies:
\[d(P, xP) \leq \mbox{vol}(\CA, P)\]

\end{lem}
\begin{proof}
(a) We only need to observe that a reduced word $w$ in $\CA$ defines a
geodesic segment $[{\bf 1}, w {\bf 1}]$ in $\cal \tilde R_{A}$, and
apply the definition of $\BBT(\CA, P)$.

\smallskip

\noindent
(b) As $w$ is cyclically reduced, $wP$ and $w^{2}P$ are contained in
the $\BBT(\CA, P)$-neighborhood of the segment $[P,w^3P]$, and thus
$P$ and $wP$ are contained in the $\BBT(\CA, P)$-neighborhood of the
segment $[w^{-1}P,w^2P]$.  As the axis $\Ax(w)$ contains $[w^{-1}P,wP]
\cap [P,w^2P]$, and the latter is the fundamental domain with respect
to the action of $w$ on $\Ax(w)$ and hence has length $\dvb w
\dvb_{T}$, the desired inequality follows.

\smallskip

\noindent
(c) As in (b) we see that $P$ and $uP$ are contained in the $\BBT(\CA,
P)$-neighborhood of the intersection $[w^{-1}P,wP] \cap [P,w^2P]$ and thus of $Ax(w)$.

\smallskip

\noindent
(d) This is a direct consequence of the above definition of the
volume.
\end{proof}

The following statement
has been shown in \cite{ll4}, Remark~2.6:

\begin{lem}
\label{lem:bbtzero}
Let $T$ be an $\R$-tree with a very small action with dense orbits of
$\FN$.  For any point $P$ in $T$, there exists a sequence of bases
${\cal A}_{k}$ such that the two sequences of constants $\BBT({\cal
A}_{k}, P)$ and $\vol({\cal A}_{k}, P)$ both tend to 0, for $k \to
\infty$.  \qed
\end{lem}


\section{The dual lamination associated to an $\R$-tree}
\label{subsec:lamtree}

In \cite{chl1-I}, \S{}2, {\em algebraic laminations} have been defined
as subsets of $\partial^{2}\FN$ (defined in \S\ref{sec:intro}
above). In particular, for every $w \in \FN$ the {\em rational}
algebraic lamination $L^{2}(w)$ has been introduced. In \S\S{}4 and 5
of \cite{chl1-I} {\em laminary languages} over a basis $\CA$ of $\FN$
have been presented, and the bijection $\rho^\CA_{\LL} \rho^{2}_{\CA}$
has been established which associates in a natural way to any
algebraic lamination $L^{2} \subset \partial^{2} \FN$ a laminary
language, denoted here by $\LL_{\CA}(L^{2})$, compare Theorem~1.1 of
\cite{chl1-I}.  We use the notation from the predecessor article
\cite{chl1-I} to facilitate the reading of this section, but we also
simplify some of the notation for the convenience of the reader. For
example, if $\LL = \LL_{\CA}(L^{2})$ is the laminary language
canonically associated to the algebraic lamination $L^{2}$, we will
write in the next section $L^{2} = L^{2}(\LL)$.

\medskip

Let $T$ be an $\R$-tree with isometric $\FN$-action.  For every
$\epsilon > 0$ we consider the set
\[
\Omega_{\epsilon}(T) \, \,  = \, \,   \{ w \in \FN
\mid  \DB w \DB_{T} < \epsilon  \} \, \,   \subset \, \, \FN
\]
which is invariant under conjugation and inversion in $\FN$, 
and
the
set
\[
\Omega^{2}_{\epsilon}(T) \, \, \, = \bigcup_{w \in
\Omega_{\epsilon}(T)} L^{2}(w) \,\,\, \subset \, \, \, \partial^{2}\FN
\]
which is invariant under the action of $\FN$ and of the flip-map on
$\partial^{2}\FN $.  We note that either the $\FN$-action on $T$ is
free simplicial, i.e. $T$ belongs to $cv_{N}$, or else
$\Omega_{\epsilon}(T)$ and hence $\Omega^{2}_{\epsilon}(T)$ are
non-empty for any $\epsilon > 0$. In the latter case we pass to the
closure in $\partial^{2}\FN $ to obtain an algebraic lamination:
\[
L^{2}_{\epsilon}(T) \,\, = \,\, \overline{\Omega^{2}_{\epsilon}(T)}
\,\, \subset \,\, \partial^{2}\FN\, .
\]
By Lemma~4.2 of \cite{chl1-I} we can define:

\begin{defn}
\label{def:lamtreeMartin}
Let $T$ be an $\R$-tree on which $\FN$ acts by isometries.  If $T$
belongs to $cv_{N}$, then we define $L^2_\Omega(T)$ to be the empty
set (which is not an algebraic lamination!).  Otherwise we define the
{\em dual algebraic lamination} associated to $T$ as follows:
\[
L_\Omega^{2}(T) = \underset{\epsilon > 0}{\bigcap} \,
L^{2}_{\epsilon}(T)\, .
\]
We note that $L_\Omega^{2}(T)$ depends only on the projective class
$[T] \in \CVNbar$.
\end{defn}

For any basis $\CA$ of $F_{N}$ we define in a similar spirit, for any
$T$ which belongs to $\partial \cvn$, the laminary languages
\[
\LL_{\CA}^{\epsilon}(T) = \bigcup_{ w \in \Omega_{\epsilon}(T)}
\LL_{\CA}(L^{2}(w))
\]
and
\[
\LL^\Omega_{\CA}(T) = \underset{\epsilon > 0}{\bigcap} \,
\LL_{\CA}^{\epsilon}(T) \, .
\]
Hence $u \in F(\CA)$ belongs to $\LL^\Omega_{\CA}(T)$ if and only if
for all $\epsilon > 0$ there exists a cyclically reduced word $w \in
F(\CA)$ with $\DB w \DB_{T} < \epsilon$ and an exponent $m \geq 1$
such that $u$ is a subword of $w^m$. However, it suffices to consider
exponents which satisfy $m \leq | u |_{\CA}\, $ (= the word length of
$u$ in $\CA^{\pm 1}$). Thus we obtain as direct consequence:

\begin{rem}
\label{mistake}
A word $u \in F(\CA)$ belongs to $\LL^\Omega_{\CA}(T)$ if and only if
for all $\epsilon > 0$ there exists a cyclically reduced word $w \in
F(\CA)$ with $\DB w \DB_{T} < \epsilon$ such that $u$ is a subword of
$w$.
\end{rem}

\smallskip

It follows easily that these laminary languages correspond precisely
to the algebraic laminations $L^{2}_{\epsilon}(T)$ and
$L_\Omega^{2}(T)$, under the bijection $\rho^\CA_{\LL} \rho^{2}_{\CA}$
established in Theorem~1.1 of \cite{chl1-I}.  Indeed, except for the
passage from $\Omega^2_{\epsilon}(T)$ to $L^{2}_{\epsilon}(T)$,
i.e. closing up in $\partial^{2}F_{N}$, the identity between the
corresponding laminary languages is definitory.  But as the language
$\LL_{\CA}^{\epsilon}(T)$ is already laminary (see \S{}5 of
\cite{chl1-I}), it follows that it agrees with the laminary language
$\rho^\CA_{\LL} \rho^{2}_{\CA}(L^{2}_{\epsilon}(T))$ associated
canonically to the closure of $\Omega^2_{\epsilon}(T)$ in
$\partial^{2} \FN$.

\bigskip

We finish this section with two observations regarding the dual 
lamination for particular cases or $\R$-trees $T$. Both of the 
following remarks follow directly from our definitions above.

\begin{rem}
\label{pointstabilizer}
Let $T \in \cvnbar$ be an $\R$-tree. For any point $x \in T$
we consider the stabilizer 
\[
\Stab(x) = \{ w \in \FN \mid wx = x \} \, .
\]
Then every $w \in \Stab(x)$ is conjugate to a word in
$\LL^\Omega_{\CA}(T)$, for any basis $\CA$ of $\FN$.

Equivalently, noting that $\Stab(x)$ is a free group $F_{k}$ of rank
$k \leq N$ (see \cite{gl}), the lamination $L^{2}_{\Omega}(T)$
contains the image (under the map canonically induced by the inclusion
$F_{k} \subset \FN$, see \cite{chl1-I}, Remark~8.1) of the full
lamination $L^{2}_{\Omega}(F_{k}) = \partial^{2}F_{k}$.
\end{rem}

The following remark is useful with respect to the canonical
decomposition of an $\R$-tree $T \in \cvnbar$ into $T_{d}$ and $T_{c}$
as given in Remark~\ref{discrete-dense}: If either of them is
non-empty, one can contract the connected components of the other one
to get $\FN$-equivariant, distance decreasing maps $T \to T/T_{d}$ or
$T \to T/T_{c}$, and both quotient $\R$-trees belong again to
$\cvnbar$.

\begin{rem}
\label{contractions}
Let $T, T' \in \cvnbar$, and let $T \to T'$ be an $\FN$-equivariant,
distance decreasing map. Then one has:
\[
L^{2}_{\Omega}(T) \subset L^{2}_{\Omega}(T')\, .
\]
\end{rem}

These two observations are the starting point for a more detailed
structural analysis of $\FN$-actions on $\R$-trees: More details will
be given in \cite{chl3}.


\section{One-sided infinite words}\label{subsec:lamtreeonesided}

Let $T$ be, as before, an $\R$-tree with a left action by isometries
of the free group $\FN$.  We fix a basis $\CA$ of $\FN$ and a point
$P$ of $T$. The choice of the basis $\cal A$ gives us an
identification between the boundary $\partial \FN$, and the space
$\partial F({\cal A})$ of (one-sided) infinite reduced words $x_{1}
x_{2} x_{3} \ldots$ in ${\cal A}^{\pm 1}$.

\smallskip

Following \cite{ll4} we denote by $L_{\cal A}^1(T) \subset \partial
F(\CA)$ the subset of those infinite reduced words $X=x_1x_2\ldots$ in
$\CA^{\pm 1}$ which have the property that for some $P \in T$ the
sequence $(X_iP)_{i\in\N}$ is bounded, where $X_i$ is the prefix of
length $i$ of $X$.  We observe:

\begin{rem}
\label{rem:monoinf}
(1) If $X$ belongs to $L_{\cal A}^1(T)$, then for {\em any} $P \in T$
the sequence $(X_{i}P)_{i\in\N}$ is bounded.

\smallskip

\noindent
(2) If $X$ does not belong to $L_{\cal A}^1(T)$, then for any $C > 0$,
any $P \in T$ and any integer $K$ there exist $l > k \geq K$ such that
$d(X_kP, X_lP) > C$.

\smallskip

\noindent
(3) If in addition $T$ satisfies BBT, then for any $X = x_{1} x_{2}
\ldots \in L_{\cal A}^1(T)$ there exists an integer $K \geq 0$ such
that for all $k, l \geq K$ one has $d(X_{k} P, X_{l}P) = d(x_{k+1}
\ldots x_{l} P, P)<3 BBT(\CA,P)$.
\end{rem}

\begin{prop}[\cite{ll4}]
\label{lun}
The canonical identification $\partial F({\cal A}) =\partial \FN$
associates to the subset $L_{\cal A}^1(T) \subset \partial F({\cal
A})$ a set $L^1(T) \subset \partial \FN$ that does not depend on the
choice of $\cal A$ (which justifies the notation $L^{1}(T)$).
\end{prop}

\begin{proof} Let $\cal B$ be another basis for $\FN$ and $X$ be in $\partial
\FN$.  Denote by $X_{\cal A}$ and $X_{\cal B}$ the corresponding
(one-sided) infinite reduced words in $\partial F({\cal A})$ and
$\partial F({\cal B})$.  The prefix sequences $(X_{{\cal{A}},
i})_{i\in\N}$ and $(X_{{\cal B}, i})_{i\in\N}$ (which are sequences of
elements of $\FN$) both converge to $X$.  Geometrically they are two
quasi-geodesics in $\FN$.  If we fix the word metric $d_{\cal A}$ on
$\FN$, then the first sequence lies on a geodesic, and the second on a
quasi-geodesic.  In particular, as follows from Cooper's cancellation
bound (see Lemma~7.1 of \cite{chl1-I}), the two sequences have the
property that for any positive integer $j$ there exists a positive
integer $i$ such that $d_{\cal A}(X_{{\cal B}, j}, X_{{\cal{A}}, i}) <
C$, where $C=\BBT({\cal B},\CA)$ is Cooper's cancellation bound
between basis $\CA$ and $\cal B$.

This shows that the sequence of points $X_{{\cal{A}}, i} P$ is of
bounded diameter if and only if the sequence of the $X_{{\cal{B}}, j}
P$ is.
\end{proof}

Let us now state some properties of $L^1(T)$: Note first that the last
proposition, and thus the definition of $L^1(T)$, does not require
that $T$ satisfies BBT.  Furthermore, the subset $L^1(T)$ of $\partial
\FN$ is $\FN$-invariant, and, unless it is empty, it is dense in
$\partial \FN$ (as is any $\FN$-invariant non-empty subset of
$\partial \FN$).  The set $L^1(T)$ is empty if and only if the action
of $\FN$ on $T$ is free and discrete.

\bigskip

Using the projection on the first coordinate, $\pi: \partial^{2}\FN
\to \partial \FN, (X, X') \mapsto X$, one can alternatively deduce
from the algebraic lamination $L^2_{\Omega}(T)$ defined in section
\ref{subsec:lamtree} a second $\FN$-invariant subset of $\partial \FN$
associated to $T$, namely the image set $\pi(L^2_{\Omega}(T))$.

\begin{prop}
\label{prop:comparison}
For any $\R$-tree $T$ with isometric $\FN$-action that satisfies BBT,
one has:
\[
\pi(L^2_\Omega(T))  \subset L^{1}(T)
\]
\end{prop}

\begin{proof}
Let $\CA$ be a basis of $\FN$ and $X = x_{1} x_{2} \ldots$ a reduced
infinite word in $\partial F(\CA)=\partial\FN$.  If $X$ does not
belong to $L^{1}(T)$, then it follows from Remark~\ref{rem:monoinf}
(2) that for any point $P\in T$, for any $C > 0$ and for any $K>0$,
there exist $k,l$ with $K\leq k\leq l$ such that $d_{\CA}(X_kP,
X_lP)>C$ where $X_i$ is the prefix of $X$ of length $i$.  It follows
directly from Lemma~\ref{lem:PBBT} (c) that, chosing $C$ large enough
(depending on the given $P$ and $\CA$), the word $X_k\inv X_l$ cannot
be a subword of any cyclically reduced word that represents an element
$w \in \Omega_{\epsilon}(T)$, for small $\epsilon > 0$.  Hence
$X_k\inv X_l$ and its inverse $X_l\inv X_k$ do not belong to
$\LL^\Omega_{\CA}(T)$, by Remark~\ref{mistake}.

But if there is some $X' = x'_{1} x'_{2} \ldots \in \partial \FN$ such
that $(X, X')$ belongs to $L^{2}_\Omega(T)$, then the reduced
biinfinite word $\rho^{2}_{\CA}(X, X') = {X}^{-1} X' = \ldots
{x}_{j+2}\inv {x}_{j+1}\inv \cdot x'_{j+1} x'_{j+2} \ldots$ in the
symbolic lamination $\rho_{\CA}^{2}(L^2_\Omega(T)) = L^\Omega_{\CA}(T)
\subset \Sigma_{\CA}$ associated to $L_{\Omega}^{2}(T)$ (compare
Proposition~4.4 of \cite{chl1-I}) will contain any of these $X_l\inv
X_k = x_{l}\inv \ldots x_{k+1}\inv$ with sufficiently large $k < l$ as
subword, in contradiction to the statement at the end of the last
paragraph.  This shows that $X$ cannot belong to
$\pi(L^2_{\Omega}(T))$.
\end{proof}

The converse inclusion with respect to
Proposition~\ref{prop:comparison} does not hold in general, as will be
seen in section~\ref{subsec:lamonsur}. In fact, one has to regard
$L^{1}(T)$ as a finer invariant of $T$ than the algebraic lamination
$L^{2}_\Omega(T)$, which only depends on a weakened version of the
topology of $T$, compare \cite{chl2}, while $L^{1}(T)$ may change when
different $\R$-tree structures are varying on a given topological tree
$T$.  For more details see \cite{chll}.  The fact that one can derive
$L_{\Omega}^2(T)$ from $L^1(T)$ will be shown below: it is a direct
consequence of Proposition~\ref{prop:recurdef}.

\medskip

For any basis $\CA$ of $\FN$ and any set $S \subset \partial \FN$ we
denote by $\LL_{\CA}^\infty (S) $ the the set of words $u \in F(\CA)$
and their inverses such that $u$ occurs infinitely often as subword in
some of the reduced infinite words $X_{\CA}$ that represent an element
$X \in S$ (we say $u$ is {\em recurrent} in $X_{\CA}$).  For any
non-empty $S$ the language $\LL_{\CA}^\infty (S)$ is laminary (because
we artificially added to $\LL_{\CA}^\infty (S)$ the inverses of any
recurrent $u$): we call it the {\em recurrent laminary language in
$\CA^{\pm 1}$ associated to} $S$.

\begin{defn}\label{def:reclam}
For any basis $\CA$ of $\FN$ and any non-empty set $S \subset \partial
\FN$ let $L^2_\infty(S)$ the algebraic lamination defined by the
recurrent laminary language in $\CA^{\pm 1}$ associated to $S$:
\[
L^2_\infty(S)=L^2(\LL_{\CA}^\infty (S))
\] 
Here (and below) we denote by $L^{2}(\LL)$ the algebraic lamination
$(\rho^\CA_{\LL} \rho^{2}_{\CA})\inv (\LL)$ defined in \cite{chl1-I}
for any laminary language $\LL$.
\end{defn}

Proposition~\ref{prop:recurdef} below justifies the absence of
mentioning explicitely the basis $\CA$ in the notation
$L^2_\infty(S)$.  But first we observe:

\medskip

For any basis $\CA$ of $\FN$, any $k \in \N$ and any $w \in F(\CA)$,
we denote by $w\chop_{k} \in F(\CA)$ the subword of $w$ obtained from
chopping off the initial and final subword of length $k$, see
\cite{chl1-I}, \S{}5.  For any second basis $\CB$ of $\FN$ consider
Cooper's cancellation bounds $C=\BBT(\CB,\CA)$ and $C'=\BBT(\CA,\CB)$
(as given in \cite{chl1-I}, \S{}7), as well as $C''=C'+\lambda C$,
where $\lambda$ is the maximal length of the elements of $\CA$ written
as words in $\CB^{\pm 1}$.

\begin{lem}
\label{rem:chop}
Consider a word $u_{\CB} \in F(\CB)$ in $\CB^{\pm 1}$, and let
$u_{\CA} \in F(\CA)$ be the word in $\CA^{\pm 1}$ representing the
same element of $\FN$. Let $v_{\CA} = u_{\CA}\chop_{C}$, and let
$w_\CA \in F(\CA)$ be any word that contains $v_\CA$ as a subword.
Let $w_{\CB} \in F(\CB)$ be the word in $\CB^{\pm 1}$ representing the
same element of $\FN$ as $w_\CA$.  Then $u_\CB\chop_{C''}$ is a
subword of $w_\CB\chop_{C'}$.
\end{lem}

\begin{proof}
This follows from a straight forward calculation.
\end{proof}

\begin{prop}
\label{prop:recurdef}
Let $S$ be a non-empty subset of $\partial \FN$, and $\CA$ and $\CB$
be two basis of $\FN$.  Then one has:
\[
L^2(\LL_{\CA}^\infty (S)) = L^2(\LL_{\CB}^\infty (S))\, .
\]
\end{prop}

\begin{proof}
It suffices to prove for any $Y \in S$ the equality
\[
L^2(\LL_{\CA}^\infty (Y))=L^2(\LL_{\CB}^\infty (Y)).
\] 
Thus, for any $(X,X')$ in $L^2(\LL_{\CA}^\infty (Y))$, we have to
prove that the biinfinite word $Z_\CB=\rho_\CB(X,X')$ in $\CB^{\pm 1}$
is contained in the symbolic lamination $L(\LL_{\CB}^\infty (Y)) =
(\rho^\CB_{\LL})\inv(\LL_{\CB}^\infty (Y))$ (see \cite{chl1-I},
\S\S{}4 and 5).  Let $u_\CB \in F(\CB)$ be a finite subword of
$Z_\CB$, and let $u_\CA \in F(\CA)$ be the word in $\CA^{\pm 1}$
representing the same element of $\FN$ as $u_\CB$. Using Cooper's
cancellation bound $C=\BBT(\CB,\CA)$, we see that
$v_\CA=u_\CA\chop_{C}$ is a subword of the biinfinite word
$Z_\CA=\rho_\CA(X,X')$.  By definition of $L^2(\LL_{\CA}^\infty (Y))$,
$v_\CA$ is recurrent in the infinite word $Y_\CA\in\partial F(\CA)$
representing $Y$.  Let $w_\CA$ be a subword of $Y_{\CA}$ such that
$v_\CA$ is a subword of $w_\CA$, sufficiently far from the beginning
and the end of $w_{\CA}$.

Let $w_{\CB} \in F(\CB)$ be the word in $\CB^{\pm 1}$ representing the
same element of $\FN$ as $w_{\CA}$.  Using Cooper's cancellation bound
$C'=\BBT(\CA,\CB)$, we get that $w_\CB\chop_{C'}$ is a subword of the
infinite word $Y_\CB\in\partial F(\CB)$ representing $Y$, and that
$u_\CB\chop_{C''}$ is a subword of $w_\CB\chop_{C'}$, for $C'' > 0 $
depending only on $\CA$ and $\CB$ as specified in
Lemma~\ref{rem:chop}.  This proves that $u_\CB\chop_{C''}$ is
recurrent in $Y_\CB$.  As $u_{\CB} \in \LL(Z_\CB)$ was chosen
arbitrarily and $\LL(Z_\CB)$ is a laminary language, this proves that
the biinfinite word $Z_\CB=\rho_\CB(X,X')$, for any $(X,X') \in
L^2(\LL_{\CA}^\infty (Y))$, is contained in the symbolic lamination
$L(\LL_{\CB}^\infty (Y))$.  Thus we have
\[
L^2(\LL_{\CA}^\infty (Y))\subset L^2(\LL_{\CB}^\infty (Y)) \, ,
\]
which proves the proposition.
\end{proof}

\begin{defn}\label{def:reclamT}
Let $T$ be an $\R$-tree with very small action of $\FN$ with dense
orbits, and let $\CA$ be any basis of $\FN$.  We define
$L^2_\infty(T)$ as the algebraic lamination defined by the recurrent
language associated to $L^1(T)$:
\[
L^2_\infty(T)=L^2(\LL_{\CA}^\infty (L^1(T)))
\] 
It follows from Proposition~\ref{lun} and
Proposition~\ref{prop:recurdef} that $L^2_\infty(T)$ does not depend
on the choice of the basis $\CA$.
\end{defn}

\begin{prop}
\label{prop:langequi}
Let $T$ be an $\R$-tree with very small action of $\FN$ with dense
orbits. The lamination $L^2_\Omega(T)$ of
Definition~\ref{def:lamtreeMartin} and the lamination $L_\infty^2(T)$
of Definition~\ref{def:reclamT} are equal:
\[
L^2_\Omega(T)=L^2_\infty(T).
\]
\end{prop}

\begin{proof}
Let $\CA$ be a basis of $\FN$. We will prove that the laminary
languages associated to these laminations (via the canonical map
$\rho^\CA_{\LL} \rho^{2}_{\CA}: \Lambda^{2} (\FN) \to
\Lambda_{\LL}(\CA)$ from \cite{chl1-I}) are equal.

\smallskip

We first prove that $\LL_{\CA}^\infty(L^{1}(T)) \subset
\LL^\Omega_{\CA}(T)$: Since the $\FN$-action on $T$ is very small and
has dense orbits, for any $\epsilon > 0$ (fixed for the rest of this
proof) there exists by Lemma~\ref{lem:bbtzero} a basis $\CB$ of $\FN$
and a point $P \in T$ such that $\BBT(\CB, P) < \epsilon$ and
$\vol(\CB,P)<\epsilon$.

Let $X = x_{1} x_{2} \ldots $ be a reduced infinite word in $\CA^{\pm
1}$ that belongs to $L^{1}_\CA(T)$.  We know that, according to the
above definition of the associated recurrent laminary language, if
$X_{k,l} = x_{k} \ldots x_{l} \in \LL_{\CA}^\infty(X)$, then there are
arbitrary large $k', l'$ with $X_{k,l} = X_{k', l'}$. Hence Cooper's
cancellation bound, see \S{}7 of \cite{chl1-I}, assures us that, when
writing $X$ as a reduced word in the basis $\CB$, say $X_{\CB} = y_{1}
y_{2} \ldots$, there is a recurrent subword $Y_{r, s} = y_{r} \ldots
y_{s}$ of $X_{\CB}$ which has the property that, when written as word
$Y_\CA$ in $\CA^{\pm 1}$, the latter contains $X_{k,l}$ as subword. On
the other hand, since $\CB$ contains at least two elements (by the
general assumption that the rank $N$ of $\FN$ is at least $2$), there
is a $y \in \CB^{\pm 1}$ such that $Y_{r, s} y$ is cyclically reduced.
Hence we obtain from Remark~\ref{rem:monoinf} (3), from
Lemma~\ref{lem:PBBT} (d) and from the above assumptions $\BBT(\CB, P)
< \epsilon$ and $\vol(\CB,P)<\epsilon$ that
\[
\DB Y_{r, s} y \DB_{T} \leq d(P, Y_{r,s} P) + d(P, yP) \leq 4
\BBT(\CB, P) < 4\epsilon\, .
\]
But for 
$s$
sufficiently large and $r$ sufficiently small  
the subword $X_{k,l}$ of
$Y_\CA$ will not be cancelled,
when $Y_{r, s} y$ is written in $\CA^{\pm 1}$
and subsequently cyclically reduced,
by Lemma~\ref{rem:chop}.
This implies that $X_{k,l}$ belongs to
$\LL_{\CA}^{4\epsilon}(T)$, which
proves the assertion.

\smallskip

We now turn to the proof of the converse inclusion, namely
$\LL^\Omega_{\CA}(T) \subset \LL_{\CA}^\infty(L^{1}(T))$: For any word
$w$ in $F(\CA)$ we distinguish between its {\em conjugating part} $v
\in F(\CA)$ and its {\em cyclically reduced part} $w' \in F(\CA)$,
where $w=vw'v\inv$ is in reduced form, with $w'$ cyclically reduced.

Let $u\in \LL^\Omega_{\CA}(T) $ be a word in $\CA^{\pm 1}$, $P$ a
point in $T$ and $\epsilon>0$.  We first want to show that there
exists a word $w\in F(\CA)$ that contains $u$ as a subword of its
cyclically reduced part and satisfies $d(P,wP)<\epsilon$.

Indeed, by Remark~\ref{mistake} there exists a cyclically reduced word
$w$ in $F(\CA)$ of which $u$ is a subword and such that $\| w\|_T<
\frac{\epsilon}{6} \, $.  Then $u$ is a subword of any cyclic
conjugate of $w^2$. As the action of $\FN$ on $T$ has dense orbits
there exists a word $v$ of $F(\CA)$ such that
\[
d( v\inv P,\Ax(w))<\frac{\epsilon}{3} \, .
\] 
An easy calculation then shows that $vw^2v\inv$ satisfies $d(P,
vw^2v\inv P) < \epsilon$, which is what we claimed.

Thus there exists a sequence of words $u_{k} \in F(\CA)$ where each of
them contains $u$ or $u\inv$ as a subword of its cyclically reduced
part, such that $d(P, u_{k} P)<\frac{1}{2^k}$.  We apply
Lemma~\ref{lem:combinatorial} stated below to obtain a sequence of
$w_{n} = u_{k_{n}}^{d_{n}}$, with $d_{n} = \pm 1$, and with the
further property that in each of the products $w_n w_{n+1}$ the
cancellation in $w_n w_{n+1}$ does not go further than the conjugating
parts of $w_n$ and $w_{n+1}$.  Then $X=w_1w_2w_3\cdots$ is a word in
$L_{\CA}^1(T)$ in which $u$ or $u\inv$ is recurrent.

This concludes the proof of Proposition \ref{prop:langequi}, modulo 
the proof of the subsequent lemma.
\end{proof}

\begin{lem}
\label{lem:combinatorial}
Let $(u_{k})_{k \in \N}$ be a sequence of words from $F(\CA)$. Then
there is an infinite subsequence $(u_{k_{n}})_{n \in \N}$ and
exponents $d_{n} = \pm 1$, such that the sequence of words $w_{n} =
u_{k_{n}}^{d_{n}}$ possesses the additional property that in each of
the products $w_nw_{n+1}$ the cancellation in $w_n w_{n+1}$ does not
exceed the conjugating parts of either $w_n$ or $w_{n+1}$.
\end{lem}

\begin{proof}
If the given words $u_{k}$ are almost all cyclically reduced, then we
can build the sequence $w_n$ by chosing inductively
$w_{n+1}={u_{n+1}}$ or $w_{n+1}={u_{n+1}}^{- 1}$ according to the
previous choice to avoid any cancellation in $w_n w_{n+1}$.

If the given sequence $(u_{k})_{k \in \N}$ contains a subsequence of
words $u_{k_{n}}$ with bounded length of their conjugating part, then
there exists a subsequence with constant conjugating part and we can
use the previous construction.

In the remaining case there exists a subsequence of words $u_{k_{n}}$
with strictly increasing length of their conjugating part.  Assume
that $w_1,\ldots,w_n$ were already chosen, $w_i=u_{k_{i}}^{d_{i}}$
with $d_{i} = \pm 1$, and with the property that the cancellation in
$w_i w_{i+1}$ is not more than the conjugating part of $w_i$ and
strictly less than the conjugating part of $w_{i+1}$.  Then replacing
the last word $w_{n}$ by its inverse $w_n\inv$ does not change this
property. If the cancellation in $w_n u_{k_{n+1}}$ is bigger than the
conjugating part of $w_n$ we replace $w_n$ by its inverse
$w_n\inv$. It follows that the cancellation in both $w_n u_{k_{n+1}}$
and $w_n u_{k_{n+1}}\inv$ is then not more than the conjugating part
of $w_n$ and strictly less than the conjugating part of $u_{k_{n+1}}$,
as the length of the latter is strictly bigger than that of the
conjugating part of $w_n$ (by our original ``strictly increasing''
condition for this case).
\end{proof}

\begin{rem}
\label{rem:monorecurrent}
The precise relationship between the various $\FN$-invariant sets of
one-sided infinite words associated to a lamination or to an $\R$-tree
is rather intricated, and it seems difficult to express the algebraic
lamination associated to an $\R$-tree properly in terms of such a set.
An attempt, however confusing or misleading it may be, is made in the
subsequent paragraph:

For any algebraic lamination $L^{2}$ we denote by
$L^{1}_{\infty}(L^{2}) \subset \partial \FN$ the set of all infinite
words $X = x_{1} x_{2} \ldots$ in some basis $\CA$ of $\FN$, which
have the property that their associated recurrent laminary language is
contained in $\LL_{\CA}(L^{2})$.  One can use Cooper's cancellation
bound (or rather a variation of Proposition~\ref{prop:recurdef}) to
show that $L^{1}_{\infty}(L^{2})$ does not depend on the choice of
$\CA$.  We believe that Proposition~\ref{prop:comparison} can be
extended to show that $\pi(L^{2})$ is always a subset of
$L^{1}_{\infty}(L^{2})$. The converse, however, seems in general to be
wrong: For example, for $\R$-tree laminations $L_{\Omega}^{2}(T)$,
where $T$ satisfies BBT and has dense orbits, we know that both of the
inclusions
\[
\pi(L_{\Omega}^{2}(T)) \subset L^{1}(T) \subset
L^{1}_{\infty}(L_{\Omega}^{2}(T))
\] 
hold, but we strongly suspect that, for any such $T$, they both are
proper inclusions.  On the other hand, we have seen above that the
three recurrent laminary languages associated to these three
$\FN$-invariant sets are equal.

\smallskip

The set $L^{1}_{\infty}(L_{\Omega}^{2}(T))$ deserves some further
attention, since it depends not on the metric on $T$, but only on the
{\em observer's topology} on $T$, compare \cite{chl2}. Contrary to
what seems to be indicated by the results of \cite{chl2}, the set
$L^{1}_{\infty}(L_{\Omega}^{2}(T))$ shows that the lamination
$L^{2}(T)$ alone suffices, without invoking the metric on $T$, to
exhibit certain completion points of the topological tree $T$ as lying
``far out at infinity''.
\end{rem}


\section{Bounded Back Tracking property and
the map $\CQ$
}
\label{subsec:bbtq}

Throughout this section we assume that $T$ satisfies the property BBT
(see \S\ref{sec:actions}), which follows for example if $T$ is small.
It is an easy exercise (compare \cite{gjll}) to show that the property
BBT ensures that every element $X \in \partial \FN \smallsetminus
L^{1} (T)$ determines, through picking any point $P \in T$ and any
sequence of elements $X_{i} \in F_{n}$ that converges towards $X$, a
well defined point
\[
\CQ(X) = \lim_{i \to \infty} X_{i} P
\]
of the Gromov boundary $\partial T$ of $T$.

\medskip

We suppose from now on that  $T$ is an $\R$-tree with very small 
$\FN$-action with dense orbits, and that for some (arbitrary) point 
$P \in T$ one has given a sequence of bases $\CA_{k}$ of $\FN$ that 
satisfies the properties assured by Lemma~\ref{lem:bbtzero}: Both 
$\vol({\cal R}_{\CA_{k}})$ and $\BBT(\CA_{k})$ tend to 0, for $k \to 
\infty$.

\medskip

For any infinite reduced word $X = x_{1} x_{2} \ldots \in \partial
F(\CA_{k})$ that represents an element of $L^{1}(T)$, the sequence of
points $(x_{1} \ldots x_{i} P)_{i\in\N}$ eventually stays in a bounded
region $R(X, k)$ of diameter $3 \BBT(\CA_{k})$ (compare
Remark~\ref{rem:monoinf} (c)), so that we can associate to $X$ a well
defined point $\CQ(X) = \lim_{k \to \infty} R(X, k)$.  It has been
shown in \cite{ll4} that $\CQ(X)$ depends only on $X \in \partial \FN$
and not on the above choice of $P$ and of the ${\cal A}_{k}$.  It is
important to note that $\CQ(X)$ may well be contained in the metric
completion $\wbar T$ of $T$, but not in $T$ itself.

\bigskip

Alternative definitions of the point $\CQ(X)$, for any $X \in \partial
\FN$, which do not need to consider an infinite change of bases of
$\FN$, have been given in \cite{ll3} and in\cite{ll4}, Lemma~3.4:

\begin{lem}
\label{lem:Qlimsup} For all $X$ in $\partial \FN$, for any sequence
of points $X_i \in \FN$ which converge towards $X$, and for any point
$P$ of $T$, one has in $\overline T \cup \partial T$:
\[
[P, \CQ(X)]=\overline{\underset{n \in \N}{\bigcup}\, \,
\underset{i\geq n} {\bigcap} \, \, [P,X_iP]}\, .
\]
\end{lem}

\begin{lem}
\label{lem:seqQ}
For all $X$ in $\partial \FN$ and for all $P$ in $T$, the point
$\CQ(X)$ is the only point of $\wbar T\cup\partial T$ such that there
exists a sequence of elements $X_i \in F_{n}$ which converge towards
$X$ and a point $P$ in $T$ such that the points $X_i P$ converge to
$\CQ(X)$.
\end{lem}

From Lemma \ref{lem:Qlimsup} and Remark~\ref{rem:monoinf} (3) it
follows directly:

\begin{lem}
\label{lem:Qtrapped}
Let $P$ be a point in an $\R$-tree $T \in \partial cv_{N}$ with dense
$\FN$-orbits, and let $\CA$ be a basis of $\FN$.  Then for every $X =
x_{1} x_{2} \ldots \in L_{\CA}^{1}(T)$ there exists a bound $K \geq 0$
such that for every $X_{k}= x_{1}\ldots x_{k}$ with $k \geq K$ one
has:
\[
\singlebox
d(X_{k} P, \CQ(X)) \, \, \leq \, \, 3 \, \BBT(\CA, P)\, .
\esinglebox
\]
\end{lem}

\bigskip

Summarizing the above discussion, we observe that for every very small
$\R$-tree $T$ with dense orbits every boundary point $X \in \partial
\FN \smallsetminus L^{1}(T)$ determines a point $\CQ(X) \in \partial
T$, while $X \in L^{1}(T)$ determines a point $\CQ(X)$ in $T$ or in
its metric completion $\wbar T$.  This defines a map
\[
{\cal Q}: \partial \FN  \to \wbar T \cup \partial T, \,  X
\mapsto \CQ(X)
\]
which is $\FN$-equivariant and surjective \cite{ll4}, but a priori it
is injective and continuous (with respect to the canonical boundary
topologies) only on $\partial F _{N} \smallsetminus L^{1}(T)$.  In
particular, on $L^{1}(T)$ the map $\cal Q$ is not continuous with
respect to the metric on $T$ (it does though possess on $L^{1}(T)$ a
kind of lower semi-continuous property, see Proposition~3.8 of
\cite{ll4}).  In \cite{chl2} it is proved that $\CQ$ is continuous if
we replace the metric topology on $T$ by the weaker {\em observer's
topology} (this is the topology for which a basis of open subsets is
given by connected components of $T \smallsetminus \{P\}$ for all
points $P$ of $T$).

The basic phenomenon for the lack of continuity of the map $\CQ$ is
illustrated as follows: If $X_{k}$ is a converging sequence of
elements from $\partial \FN$ with the property that for some point $Q
\in T$ the segments $[Q, \CQ(X_{k})]$ have pairwise intersections of
length converging to $0$, then $X = \lim X_{k}$ satisfies $Q =
\CQ(X)$, while the lengths of the segments $[Q, \CQ(X_{k})]$ may well
not converge to $0$ or even converge to $\infty$.

However, one can prove that the map $\CQ$ has the ``closed graph
property" (which will not be used in the sequel):

\begin{rem}
\label{lem:sixtwo}
Let $T$ be an $\R$-tree in $\partial\cvn$ which has dense
$\FN$-orbits, and, consider a sequence of boundary points
$X_{k}\in\partial\FN$ that converge to some $X\in\partial\FN$.  Assume
that the image points $\CQ(X_{k})\in \bar T\cup\partial T$ converge to
a point $R \in \bar T\cup\partial T$.  Then one has:
\[
R = \CQ(X).
\]
\end{rem}


\section{Geodesic lamination on a surface}
\label{subsec:lamonsur}

To gain some geometric intuition, let us consider in this section the
special case of an $\R$-tree that is dual to a measured lamination in
a surface: As in \S{}3 of \cite{chl1-I} we denote by $S$ a surface
with non-empty boundary and with negative Euler characteristic,
provided with a hyperbolic structure.  The latter is given by an
identification of the universal covering $\tilde S$ with a convex part
of the hyperbolic plane $\Hy^{2}$, which realizes the deck
transformation action of $\FN = \pi_{1}S$ on $\tilde S$ by hyperbolic
isometries. Then any geodesic lamination $\LG$ on $S$ defines, by
taking the full preimage, a geodesic lamination $\tilde \LG$ in
$\tilde S \subset \Hy^2$, on which $\FN = \pi_{1}S$ acts canonically.

There is a canonical {\em dual tree} $T_{\LG}$ with $\FN$-action by
homeomorphisms associated to $\LG$ (or to $\tilde \LG$), which is
defined by associating to every non-boundary leaf of $\tilde \LG$ a
point of $T_{\LG}$ which is not a branch point, and to the closure of
any complementary component of $\tilde \LG$ in $\tilde S$ a branch
point of $T_{\LG}$.  This association is a bijective and can be made
continuous: If the lamination $\LG$ is finite, then $T_{\LG}$ is
simplicial, so that there is no ambiguity. If $\LG$ is infinite, then
defining properly the topology of $T_{\LG}$ is much more delicate; we
refer the interested reader to \cite{chl2} where this problem has been
dealt with properly.

We now assume that the lamination $\LG$ is provided with a transverse
measure $\mu$ (see \cite{flp}).  Then the lift $\tilde \mu$ of $\mu$
to $\tilde \LG$ gives rise to a metric on $T_{\LG}$ by defining for
any points $x, y \in T_{\cal L}$, corresponding to leaves $l_{x},
l_{y} \in {\tilde \LG}$, the distance $d(x, y) = \tilde\mu(\alpha)$,
where $\alpha$ is an arc in $\tilde S$ with one endpoint on $l_{x}$
and the other on $l_{y}$, and $\alpha$ is assumed to be geodesic in
$\Hy^{2}$ and hence transverse to ${\tilde \LG}$.  This makes
$T_{\LG}$ into an $\R$-tree $T_{\mu}$ with isometric $\FN$-action.  It
is noteworthy that, in the exceptional but fascinating case where
$\LG$ is not uniquely ergodic, projectively different transverse
measures $\mu$ will produce projectively distinct $\R$-trees
$T_{\mu}$, and that the simplex of projective measures on $\cal L$
(located on the Thurston boundary of Teichm\"uller space) gives rise
to an anologous simplex of $\R$-trees in $\partial CV_{n}$.

We now consider an arbitrary point $Q$ on the boundary $\partial
\tilde S \subset S_\infty^1 = \partial \Hy^2$, where $\partial \tilde
S$ also coincides via our above identification $\FN = \pi_{1}S$ with
the Gromov boundary $\partial \FN$.  Let $\beta$ be the geodesic arc
which connects some arbitrary chosen point $P$ in $\tilde S$ to $Q$.
We distinguish three cases:

\begin{enumerate}
\item
$Q$ is the endpoint of a leaf $l$ of $\tilde \LG$.  Then $\tilde
\mu(\beta)$ is finite. In fact, $\beta$ projects to a segment in
$T_{\mu}$ of length $\tilde \mu(\beta)$.  Denote by $\widehat Q \in
T_{\mu}$ the image of $Q$ under this projection.

\item
The measure $\tilde \mu(\beta)$ is infinite. Then $\beta$ projects to
an infinite arc in $T_{\mu}$, and $Q$ defines a point $\widehat Q$ in
the Gromov boundary $\partial T_{\mu}$ of $T_{\mu}$ (which is
independent of the choice of $\beta$).

\item
In the remaining case the point $Q$ defines a point $\widehat Q$ in
the metric completion $\bar T_{\mu}$ of $T_{\mu}$, and the arc $\beta$
projects to a finite open arc in $T_{\mu}$ which becomes closed when
adding the point $\widehat Q$.
\end{enumerate}

In each of the three cases the geometrically described point $\widehat
Q$ is precisely the image $\CQ(Q)$ of the point $Q$, if $Q$ is viewed
as element of $\partial \FN$ via the above identification $\partial
\tilde S = \partial \FN$.

\medskip

We would like to note that this third class is non-empty for many
lamina\-tions $\LG$: For example, if $\LG$ is the contracting (or
expanding) lamination fixed by a pseudo-Anosov automorphism $\varphi$
of $S$, then it suffices to consider a lift $\tilde \varphi$ of
$\varphi$ to $\tilde S$ that does not fix any leaf (or permute
finitely many leaves) of ${\tilde \LG}$. The repulsive fixed point of
$\tilde \varphi$ on $\partial \tilde S$ will then define such a point
$Q$.  Note that the existence of such lifts $\tilde \varphi$ of
$\varphi$ has been proved in \cite{ll1}.

\medskip

The distinction of these three cases is illuminating in that it shows
that the set $L^{1}(T_{\mu})$, given here by the cases~1 and 3, may
well be strictly larger than the $\FN$-invariant set
$\pi(L^{2}(T_{\mu})) \subset \partial \FN$ canonically associated to
$L^{2}(T_{\mu})$, given here by case~1.  Indeed, while
$L^{2}(T_{\mu})$ (and accordingly, the occurence of the case~1 above)
depends only on $\LG$ and not on $\mu$, we do not know whether (but
suspect that) the partition of the complement into cases~2 and 3 may
actually depend on $\mu$.

\medskip

The distinction of $L^{1}(T_{\mu})$ into cases~1 and 3 was the
original motivation behind the definition of the algebraic lamination
$L^{2}(T)$ given in this paper.


\section{The lamination $L_{\CQ}(T)$ and the map $\CQ^{2}$}\label{sec:LQetQ2}

By restricting the domain and the range of the map ${\cal Q}$
introduced in section~\ref{subsec:bbtq}, one obtains the map:
\[
\CQ^{1}: L^{1}(T) \to \wbar T, \, X \mapsto \CQ(X)\, .
\]
Recall from \S\ref{subsec:bbtq} that the map $\CQ^{1}$ is surjective
(see \cite{ll4}), but in general not injective.  Unless $\CQ(X)$ has a
non-trivial stabilizer in $\FN$, the map $\CQ^{1}$ is conjectured to
be finite to one (see Remark~3.6 of \cite{ll4}).

\medskip

\begin{defn}
\label{def:qlamtree}
To every very small $\R$-tree $T$ with dense orbits we associate the
following $\FN$- and flip-invariant subset of $\partial^{2}\FN$:
\[
L^2_\CQ(T) = \{ (X, X')\in\partial^2\FN  \, \mid \, \CQ(X) = \CQ(X') \}.
\]
\end{defn}

Note that, as $\CQ$ is injective on $\partial \FN \smallsetminus
L^1(T)$, for $(X,X')$ in $L^2_\CQ(T)$ one has that $X$ and $X'$ belong
to $L^{1}(T)$.  Equivalently, we know that $\CQ(X)=\CQ(X')$ lies in
$\wbar T$ and not in $\partial T$.

\begin{defn}
\label{def:Q2}
We define a map
\[
\CQ^{2}: L^2_{\CQ}(T) \to \wbar T, \, \, (X, X') \mapsto \CQ(X) =
\CQ(X')
\]
which is $\FN$-equivariant and flip-invariant.
\end{defn}

Just as remarked in \S\ref{subsec:bbtq} for the map $\CQ: \partial F_N
\to \wbar T\cup\partial T$, the above map $\CQ^{1}: L^{1}(T) \to \wbar
T$ is in general not continuous.  The relevance of the set
$L^{2}_{\CQ}(T)$ is underlined by the following:

\begin{prop}
\label{prop:Q2continuous}
Let $T$ be an $\R$-tree with a very small action of $F_N$ with dense
orbits.  Then the subset $L^2_{\CQ}(T)$ of $\partial^{2} \FN$ is closed,
and the map $\CQ^2: L^2_{\CQ}(T)\to \wbar T$ is continuous (for the metric
topology on $\wbar T$).
\end{prop}

\begin{proof} Let $(X_n,Y_n)_{n\in \N}$ be a sequence of points
from $L^{2}_{\CQ}(T)$ that converge to $(X,Y) \in \partial^{2}\FN$.
Let $\cal A$ be a basis of $F_N$, $P$ be any point in $T$, and define
$C=3 \BBT({\cal A},P)$.

Let $x_n$ the largest common prefix of the infinite reduced words in
${\cal A}^{\pm 1}$ representing $X_n$ and $X$, $y_n$ that of $Y_n$ and
$Y$, and $h$ that of $X$ and $Y$.  From the assumption $(X,Y) \in
\partial^{2}\FN$ we know that $X$ and $Y$ are different, so that $h$
is a finite word: $h \in F(\CA)$.  The assumption $(X_n,Y_n)
\overset{n \to \infty} {\longrightarrow} (X,Y)$ implies that the $X_n$
and the $Y_n$ converge to $X$ and $Y$ respectively.  Hence, for $n$
big enough one obtains that $h$ is a prefix of both, $x_n$ and
$y_n$. Indeed, since $h$ is the longest common prefix of $X$ and $Y$,
it must also be that of $X_n$ and $Y_n$.  By hypothesis one has
$\CQ(X_n)=\CQ(Y_n)$, so that the BBT property together with
Lemma~\ref{lem:Qlimsup} ensures that $hP$ lies in a $C$-neighboorhood
of $\CQ(X_n)=\CQ(Y_n)$.  But then, by the definition of $L^{1}(T)$,
the hypothesis $X_n \overset{n \to \infty}{\longrightarrow} X$ and
$Y_n \overset{n \to \infty}{\longrightarrow} Y$ implies that $X$ and
$Y$ belong to $L^{1}(T)$.  Hence Lemma~\ref{lem:Qtrapped} shows that
$\CQ(X)$ as well as $\CQ(Y)$ are contained in a $2C$-neighborhood of
$hP$.  Hence passing over to $P$ and $\CA$ with arbitrary small
$\BBT({\cal A},P)$ proves that $\CQ(X_n)=\CQ(Y_n)$ converges to
$\CQ(X)=\CQ(Y)$.
\end{proof}

As a direct consequence of Proposition~\ref{prop:Q2continuous}
we obtain:

\begin{cor}
\label{cor:lam}
The set $L^{2}_{\CQ}(T) \subset \partial^{2}\FN$ is an algebraic
lamination.
\end{cor}

We remark that, contrary to the surface lamination case, the
lamination $L^{2}_\CQ(T)$ contains all {\em diagonal leaves}: From the
definition it follows directly that, if $(X,X')$ and $(X',X'')$ are in
$L^{2}_\CQ(T)$ and $X\neq X''$, then $(X,X'')$ is also in
$L^{2}_\CQ(T)$.

\begin{prop}
\label{prop:lamtree}
Let $T$ be an $\R$-tree with a very small $\FN$-action with dense
orbits. Then the lamination $L^2_{\Omega}(T)$ of
Definition~\ref{def:lamtreeMartin} and the lamination $L^2_{\CQ}(T)$
of Definition~\ref{def:qlamtree} are equal:
\[
L^2_{\Omega}(T) = L^2_{\CQ}(T)
\]
\end{prop}

\begin{proof}
In order to show the inclusion $L^{2}_{\CQ}(T) \subset
L^{2}_{\Omega}(T)$ it suffices, in view of Theorem~1.1 of
\cite{chl1-I}, to show that for some basis $\CA$ of $\FN$ any word $z
= x_{1} \ldots x_{s}$ in the laminary language $\LL_{\CA}^{\CQ}(T)$
associated to the algebraic lamination $L^{2}_{\CQ}(T)$ (via the
canonical map $\rho^\CA_{\LL} \rho^{2}_{\CA}: \Lambda^{2}(\FN) \to
\Lambda_{\LL}(\CA)$, see \cite{chl1-I}) is also contained in the
laminary language $\LL^\Omega_{\CA}(T)$ from \S\ref{subsec:lamtree}:
By extendability of any word in a laminary language we find, for any
$k \geq 0$, a ``superword'' $u_\CA = x_{-k} \ldots x_{1} \ldots x_{s}
\ldots x_{s+k}$ in $\LL_{\CA}^{\CQ}(T)$. Thus, for $k$ sufficiently
large, we can use Cooper's cancellation bound and pass to another
basis $\CB$ with $\BBT(\CB, P) < \epsilon$ and $vol(\CB, P) <
\epsilon$, for small $\epsilon > 0$, such that the word $u_\CB$ in the
basis $\CB$, which represents the same element of $\FN$ as $u_\CA$,
contains a subword $v_\CB$ with the following properties (compare the
similar situation considered in Lemma~\ref{rem:chop}):

On one hand $v_\CB$ belongs to the laminary language
$\LL_{\CB}^{\CQ}(T)$, which implies by Lemma \ref{lem:PBBT} (a) and
Lemma~\ref{lem:Qlimsup} that $ d(P, v_\CB P) \leq 3 \BBT(\CB, P) < 3
\epsilon$.  As the rank of our free group $\FN$ satisfies $N \geq 2$,
we find an element $y \in \CB^{\pm 1}$ such that $w_\CB = v_\CB y$ is
cyclically reduced and satisfies furthermore, by Lemma~\ref{lem:PBBT}
(d), that $\dvb w_\CB \dvb_T \leq d(P, v_\CB P) + d(P, yP) \leq 3
\BBT(\CB, P) + \vol(\CB, P) < 4 \epsilon$.

On the other hand, by a double application of Cooper's cancellation
bound we obtain, since $w_\CB = v_\CB y$ is cyclically reduced, for
sufficiently large $k$, that the word $w_\CA$ in the basis $\CA$,
which represents the same element of $\FN$ as $w_\CB$, contains the
originally chosen word $z = x_{1} \ldots x_{s} \in \LL_{\CA}^{\CQ}(T)$
as subword in its cyclically reduced part.  This shows that $z = x_{1}
\ldots x_{s}$ belongs to $\LL_{\CA}^{4\epsilon}(T)$, and hence, as
$\epsilon > 0$ was picked arbitrarily small, to $\LL^\Omega_{\CA}(T)$.

\smallskip

For the converse inclusion $L^{2}_{\Omega}(T) \subset L^{2}_{\CQ}(T)$
we consider any $(X, X') \in L^{2}_{\Omega}(T)$ and apply
Proposition~\ref{prop:comparison} to deduce that both, $X$ and $X'$
belong to $L^{1}(T)$.  Hence there are well defined points $R =
\CQ(X)$ and $Q = \CQ(X')$ in $\bar T$.  If $R \neq Q$, then for any
small $\epsilon > 0$ we can pass to a basis $\CB$ of $\FN$ with $12
\BBT(\CB, P) < 4 \epsilon < d(R, Q)$, and we consider the biinfinite
word $Z_\CB = \ldots y_{i-1} y_{i} y_{i+1} \ldots = \rho_{\CB}(X, X')$
in this basis.  Any subword $w = y_{-k} \ldots y_{l}$ of $Z_\CB$ with
$k$ and $l$ sufficiently large satisfies $d(y_{1} \ldots y_{l} P, R) <
3 \BBT(\CB, P) < \epsilon$, as well as $d((y_{-k} \ldots y_{0})^{-1}
P, Q) < 3 \BBT(\CB, P) < \epsilon$, as follows from
Lemma~\ref{lem:Qtrapped}.

As $d(P, wP) = d(y_{1} \ldots y_{l} P, (y_{-k} \ldots y_{0})^{-1} P)$,
this gives $4 \epsilon < d(R, Q) \leq d(y_{1} \ldots y_{l} P, R) +
d(P, wP) + d((y_{-k} \ldots y_{0})^{-1} P, Q) < d(P, wP) + 2
\epsilon$, and thus $d(P, wP) > 2 \epsilon$.  Hence, if $w$ is subword
of any cyclically reduced word $W$ in $\CB$, we obtain from Lemma
\ref{lem:PBBT} (c) the inequalities $\dvb W \dvb_T \geq d(P, wP) - 2
\BBT(\CB, P) \geq \epsilon$.  This shows that $w$ does not belong to
$\LL_{\CB}^\epsilon(T)$, for $\epsilon > 0$ small, and hence not to
$\LL_{\CB}(T)$, contradicting the assumption $Z_{\CB} \in
\rho^{2}_{\CB}(L^2_\Omega(T))$.  Thus we have proved that $R = Q$,
which shows that $(X, X')$ belongs to $L_{\CQ}(T)$.
\end{proof}


\section{From  the boundary of 
Outer space to the space of laminations}
\label{subsec:lambdatwo}

The group $\Out(\FN )$ acts canonically (from the left) on the space
$\Lambda^{2}(\FN)$ of algebraic laminations (see \cite{chl1-I},
\S{}8), but it also acts (from the right!)  on the space $cv_{N}$ and
on its ``Thurston boundary" $\partial cv_{N}$, and this induces an
action on $CV_{N} \cup \partial CV_{N}$ (see \S\ref{sec:actions}).
This right action is defined as follows: For any $\alpha \in \Aut
(\FN)$ and any tree $T \in \cvnbar$, the length function of the image
tree $T \alpha_{*}$ is given by
\[
\dvb w \dvb_{T\alpha_{*}} = \dvb \alpha (w) \dvb_{T} \quad {\rm for
\, \, every} \quad w \in \FN\, .
\]

\begin{prop}
\label{laminationequivariant}
The map
\[
\lambda^{2}: \partial CV_{N} \to \Lambda^{2}(\FN), \,\, [T] \mapsto
L^{2}(T)
\]
is $\Out(\FN)$-anti-equivariant: For any automorphism $\alpha$ of
$\FN$ and any $[T] \in \partial CV_{N}$ one has $\alpha\inv(L^{2}(T))
= L^{2}(T\alpha_{*})$.
\end{prop}

\begin{proof*}
This follows directly from the definition of $L^{2}(T) =
L^{2}_{\Omega}(T) = \underset{\epsilon > 0}{\cap} \,
L^{2}_{\epsilon}(T) $ in \S\ref{subsec:lamtree}, since we observe, for
every $\varepsilon > 0$, that
\[
\alpha\inv(\Omega_\varepsilon(T)) = \Omega_\varepsilon(T \alpha_{*})
\]
and hence
\[
\singlebox
\alpha\inv(L^{2}_\varepsilon(T)) = L^{2}_\varepsilon(T \alpha_{*}) \, .
\esinglebox
\]
\end{proof*}

However, it is important to point out that the map $\lambda^{2}$ is
not continuous. For example, consider the case $F_{3} = F(\{a, b, c
\})$, and let $D$ be the {\it Dehn twist} automorphism given by $a
\mapsto b a b^{-1}, b\mapsto b, c \mapsto c$. Then in the last section
of \cite{cl2} it is shown that every tree $T \in cv_{3} \cup \partial
cv_{3}$ with $\dvb b \dvb_{T} > 0$ converges projectively under
iteration of $D$ to the simplicial tree $\tilde \Gamma_{b}$ which is
the Bass-Serre tree of the graph of groups decomposition $F(a,b,c) =
\, $ $<a, b> *_{<b = b'>} < b' , c >$. Thus, if $T$ is, for example,
the simplicial tree obtained from the rose ${\cal R}_{\{a, b, c\}}$ by
contracting in the universal covering the edge labelled $c$
equivariantly (while leaving all other edges of same length), one gets
that $L^{2}(T) = L^{2}(c)$, compare Remark~\ref{pointstabilizer}.
Replacing $a$ by $b^m a b^{-m}$ defines a family of new trees $T_{m}
\in \partial cv_{3}$, which projectivize to points $[T_{m}] \in
\partial CV_{3}$ that give precisely the $D$-orbit of $[T] \in
\partial CV_{3}$, which has $[\tilde \Gamma_{b}]$ as forward and
backward limit point. On the other hand, we obtain from
Remark~\ref{pointstabilizer} that $L^{2}(T_{m}) = L^{2}(c)$ for all $m
\in \Z$.  Since $L^{2} (\tilde \Gamma_{b})$ contains $L^{2}(c)$, but
is much larger (for example $L^{2}(a)$ is equally contained in $L^{2}
(\tilde \Gamma_{b})$), this example illustrates that there exist
sequences of trees $T_{k}$ converging to $T$ such that $L^{2}(T_{k})$
converges into $L^{2}(T)$, but not to $L^{2}(T)$.

\begin{rem}
\label{chltrois}
The above example is in fact typical in that the dual lamination of
the limit tree contains but is in general bigger than the limit of the
dual laminations, for a convergent sequence of $\R$-trees from
$\CVNbar$.  More details and a precise statement of this fact is given
in \cite{chl3}.
\end{rem}

\affiliationone{
Thierry Coulbois, Arnaud Hilion and Martin Lustig\\
Math\'ematiques (LATP)\\
Universit\'e Paul C\'ezanne -- Aix-Marseille III\\
av. escadrille Normandie-Ni\'emen\\
13397 Marseille 20\\ 
France
\email{Thierry.Coulbois@univ-cezanne.fr\\
Arnaud.Hilion@univ-cezanne.fr\\
Martin.Lustig@univ-cezanne.fr\\
}
}

\end{document}